\documentclass[final]{siamltex}
\newtheorem{thm}{Theorem}[section]
\newtheorem{rem}[thm]{\bf {Remark} }
\usepackage{amsfonts}
\usepackage{amssymb}
\usepackage{titlesec}
\titleformat{\section}
  {\normalfont\Large\bfseries\centering}{\thesection}{1em}{}

\title{Hilbert schemes of rational curves on Fano hypersurfaces}

\author{ B. Wang\\
Jan 3, 2015}

\begin{document}

\maketitle

\begin{abstract} In this paper we try to further explore the linear model of the moduli of rational maps. Our attempt yields following results.
Let $X\subset \mathbf P^n$ be a generic hypersurface of degree
$h$. Let $R_d(X, h)$ denote the open set of the Hilbert scheme  parameterizing   irreducible rational curves of degree $d$ on $X$.
We obtain that
  \par
(1)  If $4\leq h\leq n-1$,  $R_d(X, h)$ is an integral, local complete intersection of dimension
\begin{equation}
(n+1-h)d+n-4.
\end{equation}

(2) If furthermore $(h^2-n)d+h\leq 0$ and  $h\geq 4$,  in addition to part (1),  $R_d(X, h)$  is also  rationally connected. 

\end{abstract}

\bigskip

\section{Introduction}

We work over the field $\mathbf C$ throughout.  
Hypersurfaces $X$ of projective space $\mathbf P^n$ can be classified into three different categories:
(1) Fano, (2) Calabi-Yau, (3) of general type. In our previous papers [11], [12], [13], we conclude that,  in all three categories,
 the normal sheaves of rational curves on general hypersurfaces have
vanishing higher cohomology groups.  This property is local. 
In this paper, we concentrate on the global properties in first category, Fano hypersurfaces.
In Fano case, we expect that the ``parameter" spaces of rational curves on $X$ has a positive dimension, and so there are plenty of 
rational curves that all have no obstruction.   First let's state the main theorem.
\bigskip

Let $\mathbf P^n$ be projective space of dimension $n$ over complex numbers $\mathbb C$.  
Let \begin{equation}
R_d(X, h)\subset \{ c: c\subset X\}
\end{equation} 
 denote the open set of the Hilbert scheme parameterizing  irreducible rational curves of degree $d$ on $X$. This is a subscheme of the Hilbert scheme of
 rational curves of degree $d$ in $\mathbf P^n$. 
\bigskip
 
\begin{theorem} (Main theorem).\par
(a)  If $4\leq h\leq n-1$,   $R_d(X, h)$ is an integral, local complete intersection of dimension
\begin{equation}
(n+1-h)d+n-4.
\end{equation}
(b) Furthermore, if $(h^2-n)d+h\leq 0$ and  $h\geq 4$, $R_d(X, h)$ is a rationally connected, integral, local complete intersection of dimension
\begin{equation}
(n+1-h)d+n-4.
\end{equation}
\end{theorem}

\bigskip

\begin{rem}  A generic member of $R_d(X, h)$ in the range of theorem 1.1 is a smooth rational curve.   
This assertion
 is not included in the theorem.  As far as we know, there are elementary proofs of the existence  for lines, but there are no elementary proof for a general degree $d$.   Clemens' proof on quintic threefolds ([2]) is  one of them  we know  for general degree $d$.     
The importance of this remark is that in the range $4\leq h\leq n-1$,  J, Harris, M. Roth and J. Starr's $R_d(X)$ is a non-empty open set of our $R_d(X, h)$ because $R_d(X, h)$ will be proved to be irreducible.  
\end{rem}

\bigskip

\subsection{Related work}\quad\par

(1) The results in our previous papers [11], [12], [13]  imply  that   $R_d(X, h)$ for $n\geq 4$ is a reduced, local complete intersection of dimension
$$
(n+1-h)d+n-4,
$$
where the negative  $(n+1-h)d+n-4$ is interpreted as the Hilbert scheme being empty, and $X$ could be a generic complete intersection. 
The main theorem in this paper only addresses the remaining parts, the  irreducibility and rational connectivity 
of $R_d(X, h)$. It is clear that the irreducibility could only occur when $h$ is relatively small.   \par
 
(2) Main theorem is an extension of results in two papers [6], [7] by J. Harris, J. Starr et al,  in which they initiated the study of the open set 
$R_d(X)$ of the Hilbert scheme parameterizing smooth rational curves of degree $d$.  This is achieved through a detailed analysis of Kontsevich's moduli spaces of stable maps. They are followed by the works of Beheshti, Kumar and others.  See section 4 for the details.

\par
(3) Our method is different from that in [6], [7].  This difference stems from the beginning choice of ``parameter" spaces of 
rational maps, i.e. a parameter space that parametrizes families of  rational maps.   They used the Kontsevich's moduli spaces of stable maps, we use a linear model of it.  Both methods analyze the 
structures of  ``parameter spaces" that extend to the ``boundaries".  
The difference is rooted in string theory's approach of fields in``non-linear model"  versus ``linear model".

\bigskip

\medskip

\subsection{Outline of the proof}\quad\smallskip

In string theory, there are two different theories, ``non linear sigma model" and ``gauged linear sigma model". 
   Kontsevich's moduli space of stable maps is a starting point of the rigorous, mathematical theory for ``non linear sigma model". 
 Our research focus on  the mathematical structures of fields in   ``gauged linear sigma model", which is also called a linear model 
of stable moduli in [4].  There is a filtration on this model which is helplessly simple on its own. However its interplay with hypersurfaces  is non trivial. The reason to use the linear model is that,  the incidence scheme of rational maps on generic hypersurface in the case of study,   is a ``mostly" smooth subscheme of a projective space. Once the scheme is smooth, everything else will follow automatically. The linear model has advantages and disadvantages when comparing with 
Kontsevich's moduli space of stable maps. Our general idea in [11], [12],  [13] and this paper is to re-organize  local coordinates of the linear model by breaking it down to finite blocks, then  to analyze each block one-by-one. This is accomplished in section 2.   This method will turn the disadvantages of the linear model to its advantages. 

\bigskip

Let $S=\mathbf P(H^0(\mathcal O_{\mathbf P^n}(h)))$ be the space of all hypersurfaces of degree $h$. 

Let
$$\mathbb C^{(n+1)(d+1)}$$ be the
vector space, 
$$( H^0(\mathcal O_{\mathbf P^1}(d))^{\oplus n+1}$$
whose open subset parametrizes  the set of maps $$\mathbf P^1\to \mathbf P^n$$ whose push-forward cycles have degree $d$.
\footnote {The automorphism of $\mathbf P^1$ induces a
$PGL(2)$ group action on $\mathbf P(\mathbb C^{(n+1)(d+1)})$.  Let $$PGL(2)( c_0)\subset \mathbf P(\mathbb C^{h(d+1)})$$ be the orbit of $c_0\in \mathbf P(\mathbb C^{h(d+1)})$. } 
Throughout the paper, we let $$M= \mathbb C^{(n+1)(d+1)}.$$ 
$M$ has affine coordinates.  The ``gauged linear sigma model" uses the space $M$ that has a stratification of closed subvarieties,
\begin{equation}
M=M_d\supset  M_{d-1}\supset \cdots\supset  M_0=\{constant\ maps\},\end{equation}
where
\begin{equation}
M_i=\{ (g c_0, \cdots, g c_{n}): g\in H^0(\mathcal O_{\mathbf P^1}(d-i)), c_j\in H^0(\mathcal O_{\mathbf P^1}(i))\}.
\end{equation}
This stratification makes it impossible to view $M$ as a space of morphisms of  the same degree $d$, i.e. $M_d\neq Hom_d(\mathbf P^1, X)$. \par

Let  $$\Gamma$$ be  the incidence scheme
\begin{equation}\{(c, f)\subset M\times S: c^\ast(f)=f(c(t))=0\}
\end{equation} 
Let $\Gamma_{f}$ be the projection of the fibre of $\Gamma$ over $f$ to $M$. 

 The natural dominant rational map, 
\begin{equation}\begin{array}{ccc}
 \Gamma_{f} &\stackrel{\mathcal R}\dashrightarrow & R_d(X, h),
\end{array}\end{equation}
reduces theorem 1.1 to showing  that $\Gamma_{f}$ is a rationally connected, integral variety of the expected dimension.
This rational map $\mathcal R$ will be constructed and verified by the results in [9], I 6.6.1, II 2.7. and in [10], prop. 0.9.   We'll discuss the details of this in section 4. 
 
Using this conversion, in the rest of the paper we concentrate on the scheme $\Gamma_f$. 
Notice $\Gamma_f$ has an induced filtration
\begin{equation}
\Gamma_f\supset ( M_{d-1}\cap \Gamma_f)\supset \cdots\supset ( M_0\cap \Gamma_f).\end{equation}

Notice by results in [9] (mentioned above)  $\mathcal R$ is regular on the inverse of 
$R_d(X, h)$ because the rational curves in $R_d(X, h)$ are all irreducible.  But it may not be regular 
on the lower stratum of (1.4). 

Then theorem 1.1 follows from the propositions on $\Gamma_f$  below . 
\bigskip

\begin{proposition}
If $4 \leq h\leq n$, then for each $d\geq 1$,  the scheme
\begin{equation}
\Gamma_{f}\backslash  M_0
\end{equation}
is smooth.
\end{proposition}
\bigskip

\begin{rem} The scheme $\Gamma_f$ is singular at the points in $M_0$.
\end{rem}
\bigskip

\begin{proposition}
If $4\leq h\leq n-1$,  then for each $d\geq 1$,  the scheme
\begin{equation}
\Gamma_{f}\backslash  M_0
\end{equation}
is connected. 

\end{proposition}

\bigskip

\begin{rem} When $h=n$ our method failed to prove the connectivity of
$\Gamma_f$. \end{rem}

\bigskip

\begin{proposition}
If $h\geq 4$ and $(h^2-n)d+h\leq 0$, then the scheme
\begin{equation}
\Gamma_{f}\backslash M_0
\end{equation}
is a rationally connected,  integral, complete intersection of $M$ defined by
\begin{equation}
f(c(t_1))=\cdots=f(c(t_{hd+1}))=0,
\end{equation}
where $t_1, \cdots, t_{hd+1}$ are any distinct points of $\mathbf P^1$. 

\end{proposition}
\bigskip

The propositions 1.3, 1.4 follow from the proposition 1.2 which follows from 
a rather plausible, but difficult lemma

\bigskip

\begin{lemma}
Let  $\mathcal G$ be the Gauss map
\begin{equation}\begin{array}{ccc}
X &\rightarrow (\mathbf P^n)^\ast.
\end{array}
\end{equation}
Let $c:\mathbf P^1\to X$ be a non-constant regular map (with an image of any degree). Assume $X$ is generic and $h\geq 4$.  Then
for generic $$(t_1, \cdots, t_h)\in Sym^{h}(\mathbf P^1), $$
$$\mathcal G (c(t_1)), \cdots, \mathcal G(c(t_h))$$
are linearly independent.

\end{lemma}

\bigskip

\bigskip

\section{Smoothness of the linear model}\quad\smallskip

Lemma 1.5 is the key to the results. Its proof lies in the heart of one difficult question that is essential to many important problems in this area.
In this paper we would not explore this difficult question, but refer it to the complete papers [11], [12], and [13]. 
Let's prove lemma 1.5. \bigskip

\begin{proof} of lemma 1.5:   We prove it by  a contradiction. Suppose there are a generic hypersurface $X_0=div(f_0)$  of degree $h$, 
 a non-constant rational map $c_0: \mathbf P^1\to X_0$, birational to its image,  and $h$ points
$$c_0(t_1), \cdots, c_0(t_h)$$
such that 
$$\mathcal G (c_0(t_1)), \cdots, \mathcal G(c_0(t_h))$$
are linearly dependent. 
Then 
\begin{equation}
dim( \mathcal G (c_0(t_1))\cap \cdots\cap \mathcal G(c_0(t_h)))\geq n-h+1
\end{equation}
and for any vector $\alpha\in \mathcal G (c_0(t_1))\cap \cdots\cap \mathcal G(c_0(t_h))$,
\begin{equation}
{\partial f_0\over \partial \alpha}|_{c_0(t)}=0,
\end{equation}
for all $t\in \mathbf P^1$.  This is because 
$$\mathcal G (c_0(t_1))\cap \cdots\cap \mathcal G(c_0(t_h))\subset \mathcal G(c_0(t))$$
for all $t$. 
Let $\{\alpha_j, j=1, \cdots,  r=n-h\}$ be a  set of linearly independent vectors in 
$$\mathcal G (c_0(t_1))\cap \cdots\cap \mathcal G(c_0(t_h))$$
Then $c_0$ lies on the hypersurfaces
\begin{equation}
{\partial f_0\over \partial \alpha_j}|_{c_0(t)}=0, j=1, \cdots, r.
\end{equation}
(Notice $f_0, {\partial f_0\over \partial \alpha_j}$ are generic in the moduli of hypersurfaces).
Hence it lies on the intersection
\begin{equation}
Y=\cap_j\{{\partial f_0\over \partial \alpha_j}=0\}\cap X_0.
\end{equation}

Next we are going to apply theorem 1.1 in [13]. We should elaborate the requirements for the theorem.
Let's denote the sequence of hypersurfaces defining the intersection $Y$ by 
\begin{equation}
f_0, f_1={\partial f_0\over \partial \alpha_1}, \cdots, f_{r}={\partial f_0\over \partial \alpha_r}.
\end{equation}
There are three requirements for the proof of theorem 1.1 of [13]:\par
(1) the subvariety defined by $f_j=0, 0\leq j\leq r$ is smooth of
dimension $n-r-1$ at $c_0(\mathbf P^1)$;\par
(2) each $f_j, j=0, \cdots, r$ is a generic hypersurface. This is different from the actual notion ``generic complete intersection" which
usually means that the point 
\begin{equation} \begin{array} {c} 
(f_1, \cdots, f_r)\in \\
H^0(\mathcal O_{\mathbf P^n}(h))\times H^0(\mathcal O_{\mathbf P^n}(h-1))\times\cdots\times H^0(\mathcal O_{\mathbf P^n}(h-1))
\end{array}\end{equation}
is generic. \par
(3) the first condition (1) says the subvariety $Y$ at $c_0(\mathbf P^1)$ is  a local complete intersection. The requirement is that
 dimension of the local complete intersection is larger than or equal to $3$.\par

First two conditions are satisfied because  $f_0$ is generic. By our assumption $h\geq 4$, we obtain that 
\begin{equation}
dim(Y)=h-1\geq 3.
\end{equation}
The thid condition is also satisfied. Therefore by the theorem 1.1 in [13], 

\begin{equation}
H^1(N_{c_0/Y})=0.
\end{equation}

Next we apply $H^1(N_{c_0/Y})=0$ to deduce an inequality. First
let
\begin{equation}
c_0^\ast(T_Y)=\mathcal O_{\mathbf P^1}(a_1)\oplus\cdots\oplus \mathcal O_{\mathbf P^1}(a_{dim(Y)}).
\end{equation}
 Because $H^1(N_{c_0/Y})=0$, \begin{equation}
a_j\geq -1, j=1, \cdots, dim(Y).
\end{equation}
Because at least one $a_j$ is larger than or equal to 2 ( from automorphisms of $\mathbf P^1$ ), we obtain that

\begin{equation}
c_1(c_0^\ast(T_Y))\geq -dim(Y)+3= -h+4.
\end{equation}

Now we use adjunction formula to find

\begin{equation}
c_1(c_0^\ast(T_Y))=[ n+1-(h+(n-h)(h-1))]d.
\end{equation}
Now we apply the inequality $h\leq n-1$ to obtain that
\begin{equation}
c_1(c_0^\ast(T_Y))\leq -h+3
\end{equation}
 Then (2.11) becomes 
\begin{equation}
-h+3\geq -h+4.
\end{equation}

This is absurd.   Therefore $c_0$ does not exist. The lemma 1.5 is proved.

\end{proof}

Next we prove proposition 1.2:\par
\begin{proof}   The idea of the proof is similar to that in [11] or [12].  We are going to choose affine coordinates for $M$ and
defining equations for $\Gamma_f$. Then use them to calculate the Jacobian matrix of $\Gamma_f$.
Let's start with coordinates of $M$. 
We consider a $c_0\in \Gamma_f\backslash  M_0$.  Let $c_0'$ be the normalization of $c_0(\mathbf P^1)$.  Then lemma 1.5 holds for 
$c_0'$.  Let $t_1, t_2, \cdots, t_h$ be those  points in lemma 1.5.,  i.e.,
$$\mathcal G (c(t_1)), \cdots, \mathcal G(c(t_h))$$
are linearly independent. Next we extend $$\mathcal G (c(t_1)), \cdots, \mathcal G(c(t_h))$$ to
coordinates of $\mathbf P^n$. That is to 
choose affine coordinates $$z_0, \cdots, z_n$$ of
$\mathbb C^{n+1}$ such that
 $\{z_i=0\}$ for $i=1, \cdots, h$ are exactly $\mathcal G (c_0(t_i))$. 
Next we choose affine coordinates for $M$. In each copy $H^0(\mathcal O_{\mathbf P^1}(d)$ of $M$, we
express $$c_j(t)=\sum_{k=0}^d c_j^k t^k \in H^0(\mathcal O_{\mathbf P^1}(d))$$ (for $j$-th copy) as

\begin{equation}
c_j(t)= \sum_{k=0}^d \theta_j^k (t-t_j)^k.
\end{equation}
where $t_j$ for $j=1, \cdots, h$ are the those in lemma 1.5, and $t_j=0$ if $j$ is not in the interval
$[1, h]$. 
The $\theta_j^k$ are affine coordinates for $M$.  We would like to use coordinates $w_j^k$ satisfying
(a linear transformation of $\theta_j^k$)
\begin{equation}\left\{\begin{array}{cc}
w_j^k=\theta_j^k,  & k\neq 0\\
w_j^0=\sum_{k=0}^d \theta_j^k (t'-t_j)^k. &
\end{array}\right.\end{equation}
where $t'$ is a generic complex number.

Let the corresponding 
coordinates for the point $c_0$ be $\bar {\theta_j^k}$. 
Next we choose defining equations of $\Gamma_f$ at $c_0$. 
Consider the following homogeneous polynomials in $W_j^k$.
\begin{equation} \left\{\begin{array}{cc}
f(c(t'))=\sum_{j=0}^n \epsilon_j W_j^0 &\\
{\partial ^j f(c(t_1))\over \partial  t^j} & j=1, \cdots, d\\
\cdots & \\
{\partial ^j f(c(t_h))\over \partial  t^j} & j=1, \cdots, d
\end{array}\right.
\end{equation} 
We claim that  these polynomials define the scheme $\Gamma_f$. 

To see this,  we let $c$ be a point in the scheme defined by the polynomials in (2.17). Also
let  
\begin{equation}
f(c(t))=\sum_{i=0}^{hd} \mathcal K_i(c) t^i.
\end{equation}
Using an automorphism of $\mathbf P^1$, we may assume $t_1=0$.
Then the equations $${\partial ^j f(c(t_1))\over \partial  t^j}=0, j=1, \cdots, d$$
imply that 
\begin{equation}
\mathcal K_i=0, 1\leq i\leq d.
\end{equation}
Then $f(c(t))$ satisfying the first set of equations
$${\partial ^j f(c(t_1))\over \partial  t^j}=0,  j=1, \cdots, d$$
becomes 
\begin{equation}
f(c(t))=\mathcal K_0(c)+r^d(\sum_{i=1}^{d(h-1)} \mathcal K_{d+i}(c) t^i).
\end{equation}
Next we repeat the same process inductively for the term 
$$\sum_{i=1}^{d(h-1)} \mathcal K_{d+i}(c) t^i$$
to obtain all $\mathcal K_i=0, i\geq 1$. At last
$\mathcal K_0=0$ because $f(c(t'))=0$.  Hence $c\in \Gamma_f$.  To prove the proposition  it suffices  to show that the Jacobian matrix of 
 $$f(c(t'))=0, 
{\partial ^j f(c(t_1))\over \partial  t^j}=0 , 
{\partial ^j f(c(t_h))\over \partial  t^j}=0 $$ with respect to the variables $w_j^k$ has full rank.
(Note $w_j^k$ are the coordinates for $c$). 
Let $\alpha_j^k$ be the variables for $T_{c_0}M$ with  the basis $ {\partial \over \partial w_j^k}$. 
Consider the subspace $V_T$ of $T_{c_0}M$ defined by  $\alpha_0^k=0=\alpha_l^k$ where   $k\neq 0$  and $l$ is not one of
$1, \cdots, h$. Then $T_{c_0}\Gamma_f\cap V_T$ consists of all $\alpha\in V_T$ satisfying
\begin{equation} \left\{\begin{array}{cc}
{\partial f(c_0(t'))\over \partial \alpha}=0  &\\
{\partial ^{j+1} f(c_0(t_1))\over \partial  t^j\partial \alpha}=0 & j=1, \cdots, d\\
\cdots & \\
{\partial ^{j+1} f(c_0(t_h))\over \partial  t^j\partial \alpha}=0 & j=1, \cdots, d
\end{array}\right.
\end{equation} 

 We start this with  $h$ equations in (2.21) in the second derivatives. They are equivalent to the equations
\begin{equation}
{\partial f(c_0(t_1))\over \partial \alpha_1}=\cdots={\partial  f(c_0(t_h))\over \partial \alpha_1}=0
\end{equation}
where $$\alpha_1\in span(\alpha_1^1, \cdots, \alpha_h^1).$$
By the lemma 1.5, we know that $${\partial f(c_0(t_i))\over \partial \alpha_j^1}=\delta_i^j$$ where 
$\delta_i^j=0$ for $i\neq j$ and $\delta_i^i\neq 0$.  Then the equations (2.22) implies
$\alpha_j^1=0, j=1, \cdots, h$.
Next step is to consider another $h$ equations in third derivatives. 
\begin{equation}
{\partial^{3} f(c_0(t_1))\over \partial t^2\partial \alpha}=\cdots={\partial ^3 f(c_0(t_h))\over \partial t^2\partial \alpha}=0
\end{equation}
Because $\alpha_j^1=0, j=0, \cdots, n$, we simplify (2.23) to
\begin{equation}
{\partial f(c_0(t_1))\over \partial \alpha_1^2}=\cdots={\partial f(c_0(t_h))\over \partial \alpha_h^2}=0.
\end{equation}

Then we use  lemma 1.5 to obtain that
\begin{equation}
\alpha_1^2=\cdots=\alpha_h^2=0.
\end{equation}
Recursively we obtain that the solution to the system of linear equations (2.21)  is
all $\alpha_j^k, j=1, \cdots, h, k=0, \cdots, d$ satisfying
\begin{equation}\begin{array}{cc}
\sum_{i=0}^n {\partial f(c_0(t'))\over \partial \alpha_j^0}=0 &\\

\alpha_j^k=0, j=1, \cdots, h, k=1, \cdots, d.\end{array}
\end{equation}
This means that the set of solutions to the equations (2.21)  has dimension 
 $h-1$.  Thus the rank of Jacobian matrix of $\Gamma_f$ at $c_0$ is $hd+1$, 
i.e. it has full rank.  Hence $\Gamma_{f}$ is smooth at
$c_0$ whenever $c_0$ is a non-constant.

This completes the proof.

\end{proof}

\bigskip

\subsection{Connectivity and a version of ``bend and break"} \quad\smallskip

This section will prove proposition 1.3.   In last section we proved that $$\Gamma_{f}\backslash M_0$$
is a smooth variety of dimension $$(n+1)(d+1)-(hd+1).$$
To show it is irreducible, it suffices to show it is connected.
The idea of the proof is to connect a generic point of $\Gamma_f$  to a point in the lower stratum. Then by the induction it is connected
to a point  parametrizing the multiple of lines. This is our version of ``bend-and-break".\footnote{ Our ``bend and break" fails when $h\geq n$. The failure  is due to a 
potential existence of certain irreducible components.  But we don't have an example of such failure. } 
Let $\Gamma_{f}'$ be an irreducible component of $\Gamma_{f}$.
Assume $d\geq 2$.

 Then
\begin{equation}
dim(\Gamma_f')=(n+1)(d+1)-(hd+1)=(n+1-h)d+n
\end{equation}
Let \begin{equation} M^{d-1}=\mathcal O_{\mathbf P^1}(d-1)^{\oplus n+1}.
\end{equation}

We should note that  $M_{d-1}\simeq \mathbb C\times M^{d-1}$.  It has a similar stratification
\begin{equation}
M^{d-1}=M_{d-1}^{d-1}\supset  M^{d-1}_{d-2}\supset \cdots\supset  M^{d-1}_0=\{constant\ maps\},\end{equation}
where
$$
M^{d-1}_i=\{ (g c_0, \cdots, g c_{n}): g\in H^0(\mathcal O_{\mathbf P^1}(d-1-i)), c_j\in H^0(\mathcal O_{\mathbf P^1}(i))\}.
$$
Then every irreducible components of $\Gamma_f'\cap M_{d-1}$ is  isomorphic  to an irreducible component
of 
\begin{equation} \mathbb C\times \Gamma^{d-1}_f, \end{equation}
where $\Gamma^{d-1}_f$ is defined to be 
\begin{equation}
\{c\in M^{d-1}: c\subset f\}, 
\end{equation}
and $\mathbb C$ is an affine open set of $\mathbf P(H^0(\mathcal O_{\mathbf P^1}(1)))$. 
Notice 
\begin{equation}\begin{array}{c}
dim(\Gamma^{d-1}_f\cap M^{d-1}_0)=d+n-1\\
dim(\Gamma^{d-1}_f)=(n+1-h)(d-1)+n
\end{array}\end{equation}
Because $h\leq n-1$, $d\geq 2$,  \begin{equation}
dim(\Gamma^{d-1}_f)>dim(\Gamma^{d-1}_f\cap M^{d-1}_0).
\end{equation}
The inequality (2.33) holds for every components of $\Gamma^{d-1}_f$. 
Therefore every component of  $\Gamma_f'\cap M_{d-1}$ contains a non-constant $c$. 
Thus  inside smooth locus of $\Gamma_f$, every point is connected to
a point in the lower stratum. Then by the induction it suffices to prove that the second lowest stratum $\Gamma_f\cap  (M_1\backslash M_0)$, which consists of all maps that correspond to lines,  is connected.  
By the classical result of Fano variety of lines,  this is correct.  More precisely  
$$\Gamma_f\cap  M_1$$
is isomorphic to
$$\mathbb C^{d-1}\times \Gamma_f^{1}$$
where $\Gamma_f^{1}$ is the same as (2.31) with $d=2$, and $\mathbb C^{d-1}$ is an affine open set of 
$$\mathbf P(H^0(\mathcal O_{\mathbf P^1}(d-1))).$$
Then it suffices to prove $$ \Gamma_f^{1}$$ is irreducible. 
The image of $ \Gamma_f^{1}$ under the rational map
$\mathcal R$ is just an open set of  Fano variety $F(X)$ of lines on the generic hypersurface $X=\{f=0\}$. 
It is connected by the classical result (see theorem 4.3, [9]). Therefore 
the proposition 1.3 is proved.\bigskip

\section{Rationally connectedness }

\begin{proof}
Note 
\begin{equation}  \mathbf P(M_0)
\end{equation}
is a smooth subvariety of $\mathbf P^{(n+1)d+n}$, with dimension
$$d+n-1.$$
Choose two generic planes $V_{top}, V_{bott}$ in $\mathbf P^{(n+1)d+n}$ with dimensions
$$nd-1, n+d$$
respectively. 
Consider the dominant projection map
\begin{equation}\begin{array}{ccc}
\Gamma_f\backslash  (\Gamma_f\cap V_{top})&\rightarrow & V_{bott}.
\end{array}\end{equation}
 
Because $d\geq 2$, the fibre's dimension  is at least $$(n-h)d-1$$ which is larger than or equal to  one.  By Bertini's theorem, the generic fibre is a smooth complete intersection
of $hd+1$ hypersurfaces of degree $h$ followed by  $n+d$ many hyperplanes in a projective space
of dimension $(n+1)d+n$. 

Notice the generic fibre satisfies
\begin{equation}
h(hd+1)+n+d\leq (n+1)d+n.
\end{equation}
(because $(h^2-n)d+h\leq 0$), where the left hand side is
the sum of the degrees of all hypersurfaces and right hand side is the dimension of the projective space. 
Hence the generic fibre is a smooth Fano variety. By V (2.1) and (2.13) of [9],  it is rationally connected.
By corollary 1.3, [5], 
\begin{equation}
\Gamma_f\backslash  (\Gamma_f\cap V_{top})
\end{equation}
must be rationally connected. 
The proof is completed.

\end{proof}

\bigskip

To summarize it, we just proved that

\bigskip

\begin{theorem} 
\quad \par
(1)  If $4\leq h\leq n-1$,   $\Gamma_f$ is an integral,  complete intersection of dimension
\begin{equation}
(n+1-h)d+n.
\end{equation}
(2) Furthermore, if $(h^2-n)d+h\leq 0$ and  $h\geq 4$, $\Gamma_f$ is a rationally connected, integral, complete intersection of dimension
\begin{equation}
(n+1-h)d+n.
\end{equation}

\end{theorem}

\bigskip

\begin{proof} of theorem 1.1.:   Next we show that the results of theorem 3.1 also hold for the open set $R_d(X, h)$ of Hilbert scheme.
Let $c_{bi}\in \Gamma_f$  be a  point of $\Gamma_f$ such that  $c_{bi}, \mathbf P^1\to X$  is birational to its image. 
There is a rational map 
map \begin{equation}\begin{array}{ccc}
\Gamma_f & \stackrel {\mathcal R_1}\dashrightarrow & Hom_{bir}(X)^{sn}\\
( c_0, \cdots, c_n) &\rightarrow & graph( \{t\}\to [c_0(t), \cdots, c_n(t)]) 
\end{array}\end{equation}
where $sn$ stands for semi-normalization and $\mathcal R$ is an local isomorphism at $c_{bi}$.   Next we use the results from [9], namely  I   theorem 6.3, II comment 2.7.  to  construct the composition
in a neighborhood of $c_{bi}$, 
\begin{equation}\begin{array} {cccc}
Hom_{bir}(X)^{sn} &\rightarrow & CH(W) &\rightarrow Hilb(X)^{sn}
\end{array}\end{equation}

Finally $\mathcal R$ is defined to be the composition in a neighborhood

\begin{equation}\begin{array} {cccccc} \Gamma_f &\stackrel{\mathcal R_1}\rightarrow &
Hom_{bir}(X)^{sn} &\stackrel{\mathcal R_2} \rightarrow & CH(W) &\stackrel{\mathcal R_3}\rightarrow Hilb(X)^{sn}.
\end{array}\end{equation}

By proposition 1.3, $c_{bi}$ is a smooth point of $\Gamma_f$. Then $Hom_{bir}(X)^{sn}$ is normal at $c_{bi}$. Then
the map $\mathcal R$ is regular at $c_{bi}$ because $\mathcal R_3$ is an isomorphism by I theorem 6.3, [9], 
$\mathcal R_1$  is
a smooth map with the fibres of dimension $1$ by the argument for (2.15), [11], and $\mathcal R_2$ is  the projection of a fibre product with fibres of dimension $3$ by II 2.7, [9], prop. 0.9 [10].   
Then theorem 1.1 follows from theorem 3.1

\end{proof}

\bigskip

\section{Work of Harris et al.}

Our main theorem extends the current known results in this area. 
Let $R_d(X)$ be the open set of the Hilbert scheme parametrizing smooth, irreducible rational curves of degree $d$.
One should notice $R_d(x)\neq R_d(X, h)$. In [6] and [7], J. Harris, J. Starr et al proved that
\bigskip

\begin{theorem} (Harris, Starr et al).\par
(1) If $h<{n+1\over 2}, n\geq 3$, 
$R_d(X)$ is an integral, locally complete intersection of the expected dimension
$$(n+1-h) d+n-4.$$
(2) If furthermore $h\leq {-1+\sqrt {4n-3}\over 2}$ and $n\geq 3$, in addition to that in part (1), $R_d(X)$ is also rationally 
connected. 
\end{theorem}
\bigskip

The part (1) which is in their first paper, is furthered by Coskun, Beheshti, Kumar and many others ([1], [3], etc).  It is conjectured by Coskun and Starr ([3)] that  
if $h\leq n\geq 4$, $R_d(X)$ is irreducible and has the expected dimension
$$(n+1-h) d+n-4.$$
\bigskip

\begin{theorem} For $h\geq 4$, theorem 1.1 recovers theorem 4.1, and furthermore  
Coskun and Starr's conjecture is correct  for $ h\leq n-1, n\geq 4$. 

\end{theorem} 
\bigskip

\begin{proof}  Kim and Pandharipande  proved the conjecture for  $h=1, 2$ ([8]). The case $h=3$ is  solved by Coskun and Starr ([3]).\footnote{ Their proof provides an evidence showing that for $h=3$,  $\Gamma_f$ is not smooth at $M_1$, i.e. the numerical condition $h\geq 4$ for proposition 1.2 is necessary.}  

The case $4\leq h\leq n-1$  follows from  theorem 1.1, part (1), because in the range of $4\leq h\leq n-1$, the variety $R_d(X)$ is strictly contained in $R_d(X, h)$. 

Because of the strict containment $R_d(X)\subsetneq R_d(X, h)$, theorem 4.1 follows from theorem 1.1.

\end{proof}

\bigskip

\begin{rem}
For $h=n$, the irreducibility of neither $R_d(X, h)$ nor $R_d(X)$ are settled.  For Calabi-Yau's case $h=n+1$, both $R_d(X)$ and $ R_d(X, h)$ could become reducible and $R_d(X)$ turns out to be a strictly contained, union of irreducible components of $ R_d(X, h)$.   But there is no evidence 
showing Calabi-Yau's situation could even partially occur in the Fano case $h=n$.  
\end{rem}

\bigskip

The following lemma deals with the rational connectedness. 
\bigskip

\begin{corollary}
\quad\par
(1) If $h\leq {-1+\sqrt{4n+1}\over 2}$ and $h\geq 4$, then for each degree $d$, 
the Hilbert scheme $R_d(X, h)$ is a rationally connected, integral, local complete intersection of the expected dimension.\par
(2) If ${-1+\sqrt{4n+1}\over 2}< h<\sqrt n$ and $ h\geq 4$, then for each degree $d\geq {h\over n-h^2}$, $R_d(X, h)$ is a rationally connected, integral, local complete intersection of the expected dimension.

\end{corollary}
\bigskip

\begin{rem} The part (1) of the corollary improves Harris and Starr's bound by a little ([7]). Part (2) reveals something new which says
$R_d(X, h)$ will not immediately become non-rationally connected as the degree of hypersurface increases. 
The range $({-1+\sqrt{4n+1}\over 2}, \sqrt n)$ for $h$ serves as a ``buffer-zone" for the rational connectedness to fadeout. 
The following is the conjectural graph of 
such a distribution 
$$\begin{array}{ccc}
RC &  {Not \ all \  RC} &  {Not \ RC }\\

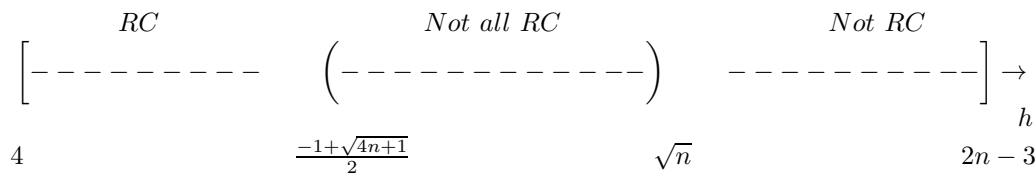
\biggl[---------&\biggl(------------\biggr)&----------\biggr]\rightarrow\\
\quad &\quad\quad  &\quad\quad\quad\quad\quad\quad\quad\quad\quad\quad\quad h \\
 4  \quad \quad \quad  \quad \quad  \quad \quad  \quad \quad   & {-1+\sqrt{4n+1}\over 2}    \quad \quad \quad  \quad \quad  \quad \quad  \quad \quad     \sqrt n & \quad \quad \quad  \quad \quad  \quad \quad  \quad \quad 2n-3
\end{array}$$
where $RC$ stands for rationally connected. In the graph
this paper proves all RC statements, but did not prove any of non RC statements for which we only know a handful of indirect examples.
For the irreducibility this ``buffer zone" may only consist of one number.
See section 5. \end{rem}
 
\bigskip

\begin{proof}
If $h\leq {-1+\sqrt{4n+1}\over 2}$, $h^2+h-n\leq 0$. Hence $h^2-n<0$. Then
\begin{equation}
(h^2-n)d+h\leq 0
\end{equation}
holds for all $d\geq 1$. Then by Main theorem 1.1, the part (1) is proved.\par

The part (2) of the corollary is just the part (2) of Main theorem 1.1.

\end{proof}

\bigskip

\section{Hilbert scheme of rational curves}
\quad

In this section we would like to organize our results in the area of rational curves on hypersurfaces. 
This extends to  hypersurfaces in other two categories: Calabi-Yau, of general type.  As before  let $X\subset \mathbf P^n$ be a generic hypersurface of degree
$h$.  Let $R_d(X, h)$ denote the open set of the Hilbert scheme  parameterizing  irreducible rational curves of degree $d$ on $X$.

In general the full scheme  structure of  $R_d(X, h)$ depends on the full scheme structure of $X$.  But we hope that  some of basic structures of
$R_d(X, h)$ may only  depend on the indices $d, h$ and $n$. We would like to discuss these structures.  
This works well in the case $h\geq 4, n\geq 5$. But as $h, n$ get smaller, the situation either becomes too simple or too complicated. They do not fit into
a bigger picture which usually has a pattern of a gradual change.   In the following we describe three basic structures of Hilbert scheme $R_d(X, h)$: (1) existence, (2) irreducibility,
(3) rational connectivity. \bigskip

(1) For the existence, we have the following demographic picture for $n\geq 5$. Its correctness  was proved in ([12]).

$$\begin{array}{ccc}
exists  &{not \ all \  exist} & {does \ not\ exist }\\
---------&\biggl [------------\biggr ]&------\stackrel {h}\rightarrow\\
\quad &\quad\quad  &\quad\quad\quad\quad\quad\quad\quad\quad\\
   & n+2    \quad \quad \quad  \quad \quad  \quad \quad  \quad \quad      2n-3 & 
\end{array}$$

 \bigskip
For $n\leq 4$, situation is subtle. If $n=4$, this is the Clemens' conjecture. We have a picture (proved in [11]), 

$$\begin{array}{ccc}
exists &   &  {does \ not\ exist }\\

---------&\biggr ]&------\stackrel {h}\rightarrow\\
\quad &\quad\quad  &\quad\quad\quad\quad\quad\quad\quad\quad\\
   & 5     & 
\end{array}$$

For $n=3$, we have a conjectural picture
$$\begin{array}{ccc}
exists & not \ all\ exist  &  {does \ not\ exist }\\

---------&\biggl [------------\biggr ]&------\stackrel {h}\rightarrow\\
\quad &\quad\quad  &\quad\quad\quad\quad\quad\quad\quad\quad\\
   & 3    \quad \quad \quad \quad \quad \quad  \quad \quad  \quad \quad  \quad \quad     4  &

\end{array}$$

(2) For the irreducibility, we have the following demographic picture for $h\geq 4, n\geq 5$, 

$$\begin{array}{ccc}
irreducible &  {not \ all \  irreducible} &  {reducible  }\\

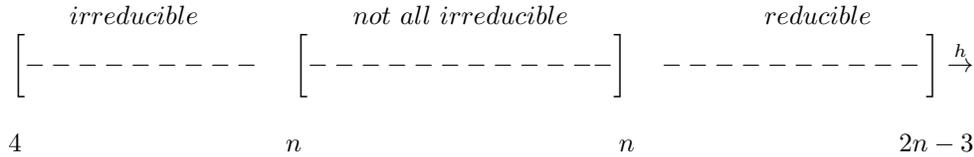
\biggl[---------&\biggl [------------\biggr]&----------\bigg]\stackrel{h} \rightarrow\\
\quad &\quad\quad  &\quad\quad\quad\quad\quad\quad\quad\quad\quad\quad\quad  \\
  4  \quad \quad \quad  \quad \quad  \quad \quad  \quad \quad  & n    \quad \quad \quad\quad \quad \quad  \quad \quad  \quad \quad  \quad \quad      n & \quad \quad \quad  \quad \quad  \quad \quad  \quad \quad 2n-3
\end{array}$$

The the statements for $h\geq n$ are our conjectures.

\bigskip

(3) For the rational connectivity, we have the conjectural demographic picture in the last section.

\bigskip

We have skipped some of the cases for  small $h, n$ (less than $4$). The
situations in those cases are subtle, and some of difficult ones are already known which do not fit 
into the larger pictures. This means for  smaller range of $n, h$ (less than $4$ ), moduli of hyperufaces 
may not be the best moduli in describing these structures.

\end{document}